\documentclass[12pt]{amsart}
\usepackage{amsmath}
\usepackage{amssymb}
\usepackage{amsfonts}
\usepackage{amsthm}
\usepackage{mathrsfs}
\usepackage{comment}
\usepackage[normalem]{ulem}
\usepackage{enumerate}
\usepackage[
textwidth=405pt,
textheight=615pt,
hmarginratio=1:1,
vmarginratio=1:1]{geometry}
\usepackage{mathtools}
\mathtoolsset{showonlyrefs}
\usepackage{hyperref}
\hypersetup{colorlinks=true,
linktoc=all,linkcolor=black,
citecolor=black}

\newcommand{\id}{\operatorname{\text{\bf I}}}

\newcommand{\hDelta}{\varDelta}

\newcommand{\C}{\mathbb{C}}
\newcommand{\D}{\mathbb{D}}
\newcommand{\T}{\mathbb{T}}
\newcommand{\Te}{\mathbb{T}}

\newcommand{\Z}{\mathbb{Z}}
\newcommand{\R}{\mathbb{R}}

\newcommand{\diff}{\mathrm{d}}
\newcommand{\Pop}{\mathbf{P}}

\newcommand{\e}{\mathrm{e}}
\newcommand{\Ordo}{\mathrm{O}}
\newcommand{\ordo}{\mathrm{o}}
\newcommand{\imag}{\mathrm{i}}

\renewcommand{\k}{\{k\}}

\newcommand{\re}{\mathrm{Re}}
\newcommand{\im}{\mathrm{Im}}
\newcommand{\Hyp}{\mathbb{H}}

\newtheorem{thm}{Theorem}[section]

\newtheorem{lem}[thm]{Lemma}

\theoremstyle{definition}

\theoremstyle{remark}

\usepackage{graphicx}

\numberwithin{equation}{section}
\title[Quasicircles and hyperbolic zero packing]
{Quasicircles and hyperbolic zero packing}
\makeatletter
\let\@wraptoccontribs\wraptoccontribs
\makeatother

\begin{document}

\thanks{This research was supported by Vetenskapsr\aa{}det
(VR grant 2020-03733), by the Leverhulme trust (grant VP1-2020-007),
and by grant 075-15-2021-602 of the Government of the Russian Federation for the 
state support of scientific research, carried out under the supervision 
of leading scientists.}

\author[Hedenmalm]{Haakan Hedenmalm}

\address{Department of Mathematics
\\
The Royal Institute of Technology
\\
S -- 100 44 Stockholm
\\
SWEDEN
\\ 
$\&$ Department  of Mathematics and Computer Sciences
\\
St. Petersburg State University
\\
St. Petersburg
\\
RUSSIA
\\
$\&$ Department of Mathematics and Statistics
\\
University of Reading
\\
Reading
\\
U.K.}
\email{haakanh@kth.se}



\date{\today}

\begin{abstract} 
We look at the work of Oleg Ivrii connected with the dimension of
quasicircles for asymptotically small quasiconformality parameter $k$.
We intend to make this work more easily accessible. Our main focus is the
integral means spectrum associated with normalized conformal mappings of
the exterior disk which have quasiconformal extensions to the whole plane
with small dilatation parameter $k$. Moreover, we address the estimates from
above only, not the sharpness from below. 
\end{abstract}

\maketitle


\section{Introduction}

\subsection{Basic notation}
\label{subsec-1.1}
We use the same basic notation as in \cite{Hed3}. So, we
write $\R$ for the real line, $\C$ for the complex plane, 
$\C_\infty:=\C\cup\{\infty\}$ for the extended complex plane 
(the Riemann sphere). We use $\diff s$ for normalized arc length
measure, and $\diff A$ for normalized area measure, while $\partial_z$,
$\bar\partial_z$ are the standard Wirtinger derivatives, and
$\hDelta_z=\partial_z\bar\partial_z$. 
We let $\D$ denote the open unit disk, $\Te:=\partial\D$ the unit circle, 
and $\D_e$ the exterior disk. We use the sesquilinear forms $\langle\cdot,
\cdot\rangle_\Te$ and $\langle \cdot,\cdot\rangle_\D$ where the measure is
$\diff s$ and $\diff A$, respectively.
We write $1_E$ for the indicator function of a subset $E$.

\subsection{The Bloch space and the Bloch seminorm}
The \emph{Bloch space} consists of those holomorphic functions 
$g:\D\to\C$ that are subject to the seminorm boundedness condition
\begin{equation}
\|g\|_{\mathcal{B}(\D)}:=\sup_{z\in\D}(1-|z|^2)|g'(z)|<+\infty.
\label{eq-Blochnorm}
\end{equation}
Let $\mathrm{aut}(\D)$ denote the group of sense-preserving M\"obius 
automorphism of $\D$. By direct calculation,
\[
\|g\circ\gamma\|_{\mathcal{B}(\D)}=\|g\|_{\mathcal{B}(\D)},\qquad
\gamma\in\mathrm{aut}(\D),
\]  
which says that the Bloch seminorm is invariant under all M\"obius 
automorphisms of $\D$.
An immediate observation we can make at this point is that provided that 
$g(0)=0$, we have the estimate
\begin{equation*}
|g(z)|\le\|g\|_{\mathcal{B}(\D)}\int_0^{|z|}\frac{\diff t}{1-t^2}
=\frac12\,\|g\|_{\mathcal{B}(\D)}\log\frac{1+|z|}{1-|z|},\qquad z\in\D,
\end{equation*}
which is sharp pointwise. This estimate is related with the interpretation
of $f\in\mathcal{B}(\D)$ as the hyperbolically Lipschitz continuous functions.

\subsection{The Bergman projection of bounded functions}
\label{subsec-PLinfty}

For $\mu\in L^1(\D)$, let  
\[
\Pop \mu(z):=\int_\D\frac{\mu(w)}{(1-z\bar w)^2}\,\diff A(w),\qquad z\in\D,
\]
be its \emph{Bergman projection}. Restricted to $L^2(\D)$, it is the 
orthogonal projection onto the subspace of holomorphic functions. In addition, 
it acts boundedly on $L^p(\D)$ for each $p$ in the interval $1<p<+\infty$ 
(see, e.g., \cite{HKZ}).
The derivative of $\Pop \mu$ is given by
\[
(\Pop \mu)'(z):=2\int_\D\frac{\bar w \mu(w)}{(1-z\bar w)^3}\,\diff A(w),
\qquad z\in\D,
\]
so that if $\mu\in L^\infty(\D)$ with $\|\mu\|_{L^\infty(\D)}\le1$, we get
that (see Per\"al\"a \cite{Per})
\begin{multline}
|(\Pop \mu)'(z)|\le 2
\int_\D\frac{|w|}{|1-z\bar w|^3}\,\diff A(w)
\\
=2
\sum_{k=0}^{+\infty}
\frac{[(3/2)_k]^2}{(\frac32+k)(k!)^2}|z|^{2k}\le 2
\Gamma(\tfrac32)^{-2}
(1-|z|^2)^{-1}=\frac{8\pi^{-1}}{1-|z|^2}.
\label{eq:Peralalaest0.0}
\end{multline}
But what happens if $m$ vanishes on a big enough disk? Can we then
get a better estimate? By the change-of-variables formula,
\[
(\Pop\mu)'(z)=2\int_\D\frac{\bar\phi(\xi)\mu\circ\phi(\xi)}
{(1-z\bar\phi(\xi))^3}|\phi'(\xi)|^2\diff A(\xi),
\]
so that with
\[
\phi(\xi)=\frac{z-\xi}{1-\bar z\xi},
\]
we find that
\[
(\Pop\mu)'(z)=\frac{2}{1-|z|^2}\int_\D\frac{\bar z-\bar\xi}
{(1-\bar z\xi)^2}\mu\circ\phi(\xi)\diff A(\xi). 
\]

\begin{lem} If $\mu\circ\phi=0$ holds on $\D(0,\varrho)$ for some $\varrho$ with
$0<\varrho<1$, then 
\begin{equation}
(1-|z|^2)|(\Pop\mu)'(z)|\le 2\int_{\D\setminus\D(0,\varrho)}\frac{1}
{|1-\bar z\xi|}\diff A(\xi)=\Ordo\bigg((1-\varrho)\log\frac{1}{1-\varrho}\bigg).
\end{equation}
\label{lem:logestimate}
\end{lem}

The details of the estimate of the integral are left to the interested reader.

\subsection{Exponential integrability}

Suppose $g=\Pop\mu$, where $\mu\in L^\infty(\D)$ with $\|\mu\|_{L^\infty(\D)}\le1$.
Let $g_r$ denotes the dilate $g_r(\zeta):=g(r\zeta)$.
It is of importance to control the integrals
\[
\int_\Te |\e^{t g_r}|\diff s
\]  
for $t\in\C$ which are close to $0$. The following estimate was offered by Ivrii \cite{Ivrii}.

\begin{thm} {\rm (Ivrii)}
In the above setting, we have the estimate
\[
\int_\Te |\e^{t g_r}|\diff s=\Ordo\bigg(\frac{1}{(1-r^2)^{(\Sigma^2+\ordo(1))|t|^2/4}}\bigg)
\]  
uniformly as $|t|\to0$. Here, $\Sigma^2$ is an absolute constant with $0<\Sigma^2\le1$.
\label{thm:Ivr}
\end{thm}

The constant $\Sigma^2$ was introduced by Astala, Ivrii, Per\"al\"a, and Prause \cite{AIPP}.
Later, Hedenmalm \cite{Hed3} found another interpretation  of the constant in terms of a zero
packing problem.

\begin{thm} {\rm (Hedenmalm)}
It holds that 
\[
\Sigma^2=1-\rho_{\mathbb{H}}<1,
\]
where
\begin{equation}
\rho_{\mathbb{H}}:=
\liminf_{r\to1^-}\inf_f\frac{\int_{\D(0,r)}\Phi_f(z)\frac{\diff A(z)}{1-|z|^2}}
{\int_{\D(0,r)}\frac{\diff A(z)}{1-|z|^2}}=
\liminf_{r\to1^-}\inf_f\frac{\int_{\D(0,r)}\Phi_f(z)\frac{\diff A(z)}{1-|z|^2}}
{\log\frac{1}{1-r^2}}.
\end{equation}
Here, the infimum runs over all polynomials $f$, and we use the notation
\[
\Phi_f(z):=\big((1-|z|^2)|f(z)|-1\big)^2,\qquad z\in\D.
\]
\label{thm:Hed}
\end{thm}

Ivrii's arXiv preprint \cite{Ivrii} has unfortunately not been published. However, this
does not mean that Ivrii's contribution is flawed. Here we shed light the various aspects
of Ivrii's argument which we believe are correct or at least may be fixed so that they are
correct.

\section{Exponential integrability}
\subsection{The Littlewood-Paley identity}

We use a variant of
the Littlewood-Paley identity, for $f\in H^2$, where $H^2$ is the usual Hardy space of
Taylor series with square summable coefficients:
\begin{multline}
\int_\Te |f|^2\diff s=|f(0)|^2+\int_\D|f'(z)|^2\log\frac{1}{|z|^2}
\diff A(z)
\\
=\int_\D|f(z)|^2\diff A(z)+\int_\D|f'(z)|^2(1-|z|^2)\diff A(z).
\label{eq:LittlewoodPaley}
\end{multline}
We apply this identity to $f(z)=\exp(\frac12 tg_r(z))$, where $g_r(z)=g(rz)$ is shorthand
for the dilated function, with $0<r<1$. 
\begin{multline}
\label{eq:LP2}
\int_\Te |\e^{t g_r}|\diff s=\int_\D|\e^{t g_r}|\diff A+
\frac{|t|^2}{4}\int_\D(1-|z|^2)|(g_r)'|^2|\e^{t g_r}|
\,\diff A(z)
\\
\le
r^{-2}\int_{\D(0,r)}\big|\e^{t g}\big|\diff A+\frac{|t|^2}{4}
\int_{\D(0,r)}(1-|z|^2)^2|g'|^2|\e^{t g}|\,
\frac{\diff A(z)}{1-|z|^2}.
\end{multline}
We may estimate the first integral on the right-hand side uniformly:
\begin{multline*}
r^{-2}\int_{\D(0,r)}|\e^{t g}|\diff A\le
\int_{\D}|\e^{t g}|\diff A\le \int_\D
\bigg(\frac{1+|z|}{1-|z|}\bigg)^{8\pi^{-1}|t|}\diff A(z)
\\
=\int_0^1
\bigg(\frac{1+r}{1-r}\bigg)^{8\pi^{-1}|t|}2r\diff r\le
\frac{4^{8\pi^{-1}|t|}}{1-8\pi^{-1}|t|},
\end{multline*}
provided that $|t|<\pi/8$. In particular, we find that
\begin{equation*}
r^{-2}\int_{\D(0,r)}\big|\e^{t g}\big|\diff A\le
\int_{\D}\big|\e^{t g}\big|\diff A\le 4,\qquad |t|\le\frac{\pi}{16}.
\end{equation*}
We implement this estimate into \eqref{eq:LP2}, and obtain that
\begin{equation}
\label{eq:LP3}
\int_\Te \big|\e^{t g_r}\big|\diff s\le
4+\frac{|t|^2}{4}
\int_{\D(0,r)}(1-|z|^2)^2|g'|^2\big|\e^{t g}\big|\,
\frac{\diff A(z)}{1-|z|^2},\qquad |t|\le\frac{\pi}{16}.
\end{equation}
By Per\"al\"a's estimate \eqref{eq:Peralalaest0.0}, the function
\begin{equation}
N_g(z):=(1-|z|^2)^2|g'(z)|^2
\label{eq:functionN_gdef}
\end{equation}
is uniformly bounded,
\begin{equation}
0\le N_g\le \frac{64}{\pi^2},
\label{eq:Peralaest0.1}
\end{equation}
and the estimate \eqref{eq:LP3} simplifies:
\begin{equation}
\label{eq:LP4}
\int_\Te |\e^{t g_r}|\diff s\le
4+\frac{|t|^2}{4}
\int_{\D(0,r)}N_g|\e^{t g}|\,
\frac{\diff A(z)}{1-|z|^2},\qquad |t|\le\frac{\pi}{16}.
\end{equation}
We proceed to analyze this inequality.

\subsection{A toy calculation}
Let us suppose for the moment that $N_g\le\eta$ holds for some positive
$\eta$. This is of course true for $\eta=64\pi^{-2}$, but we are now
interested in smaller values of $\eta$. Under this hypothesis, we see from
\eqref{eq:LP4} that
\begin{equation}
\label{eq:LP5}
\int_\Te |\e^{t g_r}|\diff s\le
4+\eta\frac{|t|^2}{4}
\int_{\D(0,r)}|\e^{t g}|\,
\frac{\diff A(z)}{1-|z|^2},\qquad |t|\le\frac{\pi}{16}.
\end{equation}
Let us write
\[
I(r,t):=\langle|\e^{tg}|\rangle_{\Te(0,r)}:=\int_\Te |\e^{t g_r}|\diff s,
\]
and observe that
\begin{equation}
\label{eq:LP6}
\int_{\D(0,r)}|\e^{t g}|\,
\frac{\diff A(z)}{1-|z|^2}=\int_0^r \frac{2u\diff u}{1-u^2}
\int_\Te|\e^{t g_{u}}|\diff s=\int_0^r I(u,t)\,\frac{2u\diff u}{1-u^2},
\end{equation}
so that \eqref{eq:LP5} would assert that for $|t|\le\frac{\pi}{16}$, 
\begin{equation}
\label{eq:LP7}
I(r,t)\le
4+\eta\frac{|t|^2}{4}
\int_0^r I(u,t)\,\frac{2u\diff u}{1-u^2}.
\end{equation}
It is possible to argue (by iteration, or in terms of differential
inequalities) that the estimate \eqref{eq:LP7} leads to the upper bound  
\begin{equation}
\label{eq:Iest0}
I(r,t)\le 4\,(1-r^2)^{-\eta|t|^2/4}.
\end{equation}
If this is true with some fixed $\eta<1$, say $\eta=\frac12$, then we would
get the desired growth estimate of $I(r,t)$ automatically. To be able to
defeat the estimate $N_g\le\eta$ the function $g'$ should exhibit some
growth. If we instead want $N_g\ge\eta$ to hold on a reasonably fat set,
the natural way to achieve the growth is by having zeros in $g'$. 
We now discuss the rather opposite behavior of the functions $\e^{tg}$
and $N_g$, respectively. Here, $\e^{tg}$ is basically hyperbolically 
locally constant, while $N_g$ is instead oscillatory.  

\subsection{The scale associated with $\e^{tg}$}
The function $g$ is in the Bloch space, and hence Lipschitz continuous
in the hyperbolic metric. Then $tg$ has a Lipschitz constant which
is $\Ordo(|t|)$, which makes it asymptotically constant on hyperbolic disks
of any fixed radius (in fact, it is basically constant over hyperbolic
distances that are $\ordo(|t|^{-1})$). On the other hand, the function $N_g$
can be much more oscillatory. 

Let us discretize the argument in the toy calculation, but without making
any fixed assertion regarding the size of $N_g$. We let $L$ be a large
positive real, which we think of as fixed, and consider successive radii
\[
r_k=\frac{1-\e^{-Lk}}{1+\e^{-Lk}},\qquad k=0,1,2,\ldots,
\]
which tend to $1$ quickly. We consider the annuli
\[
\mathbb{A}_k:=\D(0,r_k)\setminus\bar\D(0,r_{k-1}),\qquad k=1,2,3,\ldots.
\]
The hyperbolic width of the annulus $\mathbb{A}_k$ equals
\[
\int_{r_{k-1}}^{r_k}\frac{2\diff u}{1-u^2}=\log\frac{1+r_k}{1-r_{k}}-
\log\frac{1+r_{k-1}}{1-r_{k-1}}=L,  
\]
while the weighted area is
\[
|\mathbb{A}_k|_{A_{-1}}:=\int_{\mathbb{A}_k}\frac{\diff A(z)}{1-|z|^2}=\log
\frac{1-r_{k-1}^2}{1-r_k^2}=L+\Ordo(\e^{-L(k-1)}).
\]
We introduce the notation
\[
\langle h\rangle_{\mathbb{A}_k}=
\frac{\int_{\mathbb{A}_k}h\frac{\diff A}{1-|z|^2}}{\int_{\mathbb{A}_k}
\frac{\diff A}{1-|z|^2}}=\frac{1}{|\mathbb{A}_k|_{A_{-1}}}
\int_{\mathbb{A}_k}h(z)\frac{\diff A(z)}{1-|z|^2}
\]
for the weighted average of $h$ on the annulus $\mathbb{A}_k$. The object that
we will focus on are the annular averages
\[
\langle|\e^{tg}|\rangle_{\mathbb{A}_k},
\]
for $t$ close to $0$. This quantity models $I(r_k,t)$. Indeed, since $\e^{tg}$
is essentially constant in the radial direction inside $\mathbb{A}_k$ for small
$|t|$, the integral is mainly over the angular direction, which gives 
that
\[
\langle|\e^{tg}|\rangle_{\mathbb{A}_k}=\e^{\Ordo(L|t|)}
I(r_k,t)=\e^{\Ordo(L|t|)}
\langle|\e^{tg}|\rangle_{\Te(0,r_k)}.
\]

\subsection{The scale associated with $N_g$}
The function $N_g$ is bounded, but it can be highly oscillatory. It varies
at much shorter distances. We are interested in maximizing the average quantity
\[
\langle N_g |\e^{tg}|\rangle_{\mathbb{A}_k}.
\]
Since $N_g$ and $|\e^{tg}|$ vary at different scales, it is tempting to think
like a probabilist 
and suggest that we have independent stochastic variables, so that 
\[
\langle N_g |\e^{tg}|\rangle_{\mathbb{A}_k}\sim
\langle N_g \rangle_{\mathbb{A}_k}\langle|\e^{tg}|\rangle_{\mathbb{A}_k}.
\]
Let us carry on and see where this line of thinking would lead us.  
We obtain from \eqref{eq:LP4} and the assumed independence that
\begin{multline*}
\e^{\Ordo(L|t|)}\langle|\e^{tg}|\rangle_{\mathbb{A}_k}\le 4+\frac{|t|^2}{4}
\int_{\D(0,r_k)}N_g|\e^{tg}|\frac{\diff A}{1-|z|^2}
\\
\sim
4+\frac{|t|^2}{4}\sum_{j=1}^{k}
|\mathbb{A}_j|_{A_{-1}}\langle N_g\rangle_{\mathbb{A}_j}\langle |\e^{tg}|
\rangle_{\mathbb{A}_j}.
\end{multline*}
Moreover, let us suppose that for big $k$,
\begin{equation}
\langle N_g\rangle_{\mathbb{A}_j}\le \sigma_g(L)^2,
\label{eq:asymptvar99}
\end{equation}
where $\sigma_g(L)^2=\sigma_g^2+\ordo(1)$ as $L\to+\infty$.  The nonnegative
parameter $\sigma_g^2$ would express how big $N_g$ can be, on average, in
a rather precise sense. Then
we would have that
\begin{multline*}
\e^{\Ordo(L|t|)}\langle|\e^{tg}|\rangle_{\mathbb{A}_k}\lesssim
4+\frac{|t|^2}{4}\sum_{j=1}^{k}
|\mathbb{A}_j|_{A_{-1}}\langle N_g\rangle_{\mathbb{A}_j}\langle |\e^{tg}|
\rangle_{\mathbb{A}_j}
\\
\le
M+\frac{|t|^2}{4}\sigma_g(L)^2\sum_{j=1}^k|\mathbb{A}_j|_{A_{-1}}
\langle |\e^{tg}|\rangle_{\mathbb{A}_j}.
\end{multline*}
Moreover, iteration would give an estimate of the type
\begin{multline}
\langle|\e^{tg}|\rangle_{\mathbb{A}_k}=\Ordo\bigg(\exp\bigg(
\e^{\Ordo(L|t|)}\frac{|t|^2\sigma_g(L)^2Lk}{4}\bigg)\bigg)
\\
=\Ordo\bigg(\frac{1}{(1-r_k^2)^{\e^{\Ordo(L|t|)}|t|^2\sigma_g(L)^2/4}}\bigg).
\label{eq:iteration0.99}
\end{multline}
This estimate along the radii $r_k$ may be extended (with minor losses) 
to all radii $r$ if we take into account the hyperbolic Lipschitz property
with small Lipschitz constant.

\section{Grid of boxes}

\subsection{Boxes of hyperbolic size $\asymp L$}
It is actually not so difficult to make the above reasoning more rigorous.
We split the annulus $\mathbb{A}_k$ into into equal ``boxes'', by cutting
in the radial direction. Each box $Q_{k,l}$ is of the type
$\theta_{l-1}\le\theta\le\theta_l$
and $r_{k-1}\le r\le r_{k}$, so that $\theta_l-\theta_{l-1}$ is independent of
$l$. We will need each box to be of the angular width $\Ordo(L\e^{-Lk})$, so
there should be on the order of magnitude $\asymp L^{-1}\e^{Lk}$ boxes.
In each box $Q=Q_{k,l}$ we would try to maximize the average
$\langle N_g\rangle_Q$. In each such box $Q$, $|e^{tg}|$ is roughly constant,
so that
\begin{equation}
\langle N_g|\e^{tg}|\rangle_Q\sim \langle N_g\rangle_{Q}\langle |e^{tg}|\rangle_Q
\label{eq:indepSV}
\end{equation}
holds. 
All we need now is to have an estimate of the form
\begin{equation}
\langle N_g\rangle_{Q}\le \sigma_g(L)^2,
\label{eq:localN_g}
\end{equation}
with $\sigma_g(L)^2\le \Sigma^2+\ordo(1)$, where $\Sigma^2<1$ is independent of
$g$. Indeed,  then we obtain from \eqref{eq:indepSV} that
\begin{multline}
\langle N_g|\e^{tg}|\rangle_{\mathbb{A}_k}=\frac{1}{\#(Q)}
\sum_{Q}\langle N_g|\e^{tg}|\rangle_{Q}
\\
\le \e^{\Ordo(L|t|)}\sigma_g(L)^2\sum_Q
\langle N_g|\e^{tg}|\rangle_{Q}=\e^{\Ordo(L|t|)}\sigma_g(L)^2
\langle|\e^{tg}|\rangle_{\mathbb{A}_k}.
\end{multline}
This would permit us to apply the iteration argument which leads to
\eqref{eq:iteration0.99}, as desired. But what kind of argument would be
able to give us \eqref{eq:localN_g}?

\subsection{Localization of $\mu$ to boxes}
Here, we recall the formula \eqref{eq:functionN_gdef}, and split
$\mu=\mu_Q+\mu_{Q^c}$, where $\mu_Q=1_Q\mu$ and $\mu_{Q^c}=1_{\D\setminus Q}\mu$.
Then $g=\Pop\mu=g_Q+g_{Q^c}$, where $g_Q=\Pop\mu_Q$ and $g_{Q^c}=\Pop\mu_{Q^c}$.
We estimate
\begin{equation*}
(1-|z|^2)|(\Pop \mu_{Q^c})'(z)|\le 2(1-|z|^2)
\int_{\D\setminus Q}\frac{|w|}{|1-z\bar w|^3}\,\diff A(w)=\Ordo(L\,\e^{-L/2})
\end{equation*}
for $z\in\frac12Q$, if we invoke Lemma \ref{lem:logestimate}, assuming that the
corresponding hyperbolic distance between $\frac12Q$ and $\D\setminus Q$ is at least $L/2$.
This means that for the size of $N_g$ in $\frac12Q$, only $\mu_Q$ matters.
Here, it is actually better for us to work with bigger subboxes $(1-\epsilon)Q$ in place
of $\frac12Q$, which we can if $L$ is assumed big enough so that $\epsilon L$ is also big.
Here, we think of $(1-\epsilon)Q$ as a correspondingly smaller box with the same mass 
center point and all other sizes shrunk by the factor $1-\epsilon$. This procedure of going from the
box $Q$ to the slighly smaller box $(1-\epsilon)Q$ is explained in some detail
in the subsection below.
The relevance of $\epsilon L$ is that it should measure the
hyperbolic distance between $(1-\epsilon)Q$ and $\D\setminus Q$, and we need $\epsilon L$ big to
be able to ensure localization in terms of Lemma \ref{lem:logestimate}.

We aim to show that $\langle N_g\rangle_{(1-\epsilon)Q}\le \sigma_g(L)^2$ holds, with
$\sigma_g(L)^2=\sigma_g^2+\ordo(1)$, and $\sigma_g^2\le\Sigma^2$.
As a first step, we find the box $Q=Q_{k,l}$ for fixed $k$ such that
$\langle N_g\rangle_{(1-\epsilon)Q}$ is maximal (we maximize over $l$), say
$\langle N_g\rangle_{(1-\epsilon)Q}=\alpha$.
Then for that $Q=Q_{\text{max}}$, $\mu_Q$ must a rather optimal, and we may
use a periodization procedure to replace the original $\mu$ in the entire
annulus by copies of $\mu_{Q_\text{max}}$ so as to make
$\langle N_g\rangle_{(1-\epsilon)Q}\approx\alpha$ for all the other boxes
in the annulus. This would make
$\langle N_g\rangle_{\mathbb{A}_k}\approx\alpha+\Ordo(\epsilon)$.
We would now like to work with optimization along all the annuli at the
same time. In fact, is better to have scale invariance as well, which
is only approximately true for the disk. The point is that if we have found
an optimal box $Q=Q_\textrm{max}$, we would like to use $\mu_Q$ to replace
$\mu$ along all the other annuli as well, by just scaling the optimal
$\mu_Q$ (like we did for rotations before). 
Neglecting the technicalities, we should then have
$\langle N_g\rangle_{\mathbb{A}_j}\approx \alpha+\Ordo(\epsilon)$ for all
big $j$. We should think of the correct value of $\alpha$ to appear in the limit
as we zoom in at the circle $\T$, when the circle begins to look like a straight line.

\subsection{Summing over boxes (or annuli)}
It remains to find out how big can be this optimal value $\alpha$.
By summing over successive annuli we obtain
\[
\int_{\D(0,r_k)}|g'(z)|^2(1-|z|^2)\diff A=\sum_{j=1}^{k}|\mathbb{A}_j|_{A_{-1}}
\langle N_g\rangle_{\mathbb{A}_j}\approx(\alpha+\Ordo(\epsilon))(Lk+\Ordo(1)),
\]
so that
\[
\frac{\int_{\D(0,r_k)}|g'(z)|^2(1-|z|^2)\diff A}
{\int_{\D(0,r_k)}(1-|z|^2)^{-1}\diff A(z)}\approx\alpha+\Ordo(\epsilon).
\]
The left-hand side expression equals the
average of $N_g$ on the disk $\D(0,r_k)$. That average would be trivially
dominated by the maximal local average $\langle N_g\rangle_Q$, but our argument
(which needs to be made more technically sophisticated) gives that the
local maximal average over big boxes is the same as the maximal average over
the disks $\D(0,r_k)$. 
By the Littlewood-Paley identity \eqref{eq:LittlewoodPaley},
\[
\int_{\D(0,r)}|g'(z)|^2(1-|z|^2)\diff A=r\int_{\T(0,r)}|g|^2\diff s
-\int_{\D(0,r)}|g|^2\diff A
\]
it follows that the optimal value of $\alpha$ equals
McMullen's \emph{asymptotic variance}: 
\[
\sigma_g^2:=\limsup_{r\to1^-}\frac{\int_\Te|g_r|^2\diff s}{\log\frac{1}{1-r^2}}.
\]
Moreover, by Theorem \ref{thm:Hed},
\[
\sup_{g}\sigma_g^2=\Sigma^2=1-\rho_{\mathbb{H}}<1,
\]
where the supremum runs over all $g=\Pop\mu$ with $\|\mu\|_{L^\infty(\D)}$.

\subsection{Wrapping up the proof} 
We return to the various boxes $Q$. We need to keep $L$ big, $\epsilon$ small, but
$\epsilon L$ big too. Then the averages $\langle N_g\rangle_{(1-\epsilon)Q}$
approximate well $\langle N_g\rangle_Q$, given that $N_g$ enjoys the bound
\eqref{eq:Peralaest0.1}. 

The careful analysis of the boxes would require a better way to represent the disk which
is stable under dilatation, not just rotation. We may model $\D$ in terms of 
logarithmic coordinates. We write $z=\e^{\imag\pi\zeta}$, where $\im\, \zeta>0$ and
$\re\, \zeta\in\R/2\Z$. The boxes are now just rectangles that form a grid. The Bergman 
projection may be expressed in these coordinates, if  $\mathcal{D}$ is the 
fundamental domain in the upper half-plane with
real part between $-1$ and $1$:
\[
\Pop \mu (\e^{\imag\pi\zeta})=\int_{\mathcal{D}}
\frac{1}{(1-\e^{\imag\pi(\zeta-\bar\xi)})^2}\mu(\e^{\imag\pi\xi})\,\e^{-2\pi\im \xi}\diff A(\xi).
\]
Similarly, we find that
\begin{multline}
(1-|z|^2)(\Pop\mu)'(z)\big|_{z:=\e^{\imag\pi\zeta}}
\\=
2(1-\e^{-2\pi\im\zeta})\int_{\mathcal{D}}
\frac{\e^{-\imag\pi\bar\xi}}{(1-\e^{\imag\pi(\zeta-\bar\xi)})^3}\,
\mu(\e^{\imag\pi\xi})\,\e^{-2\pi\im \xi}\diff A(\xi),
\label{eq:Ngbasic}
\end{multline}
and we observe that for $\zeta,\xi\in\mathcal{D}$, the formula
\begin{multline}
\label{eq:singres0}
\frac{1}{(\e^{\imag\pi(\zeta-\bar\xi)}-1)^3}=
\\
\sum_{j=-1}^1\bigg\{(\imag\pi(\zeta-\bar\xi+2j))^{-3}-
\frac32((\imag\pi(\zeta-\bar\xi+2j))^{-2}+(\imag\pi(\zeta-\bar\xi+2j))^{-1}\bigg\}+\Ordo(1)
\end{multline}
allows us to control the singularity. The sum over $j\in\{-1,0,1\}$ is due to the fact that the
left-hand side expression is $2$-periodic in the variable $\zeta-\bar\xi$. But we may as well 
focus our attention to the
contribution with $j=0$, and observe that the contributions associated with 
the second a third term on the right-hand side of 
\eqref{eq:singres0} are negligible, i.e. correspond to asymptotically vanishing contributions
in \eqref{eq:Ngbasic} as $\im\,\zeta\to0$ (i.e., as $|z|\to1$). Similarly,
\[
\e^{-\imag\pi\bar\xi}\,\e^{-2\pi\im\, \xi}=1-\imag\pi\bar\xi -2\pi\im\,\xi+\Ordo(|\xi|^2),
\]
where only the constant $1$ makes a significant contribution. For $\zeta\in\mathcal{D}$ with
$|\zeta|\le\frac12$, we see from
\eqref{eq:Ngbasic} that
\begin{multline}
(1-|z|^2)(\Pop\mu)'(z)\big|_{z:=\e^{\imag\pi\zeta}}
\\=
-4\pi(\im\,\zeta+\Ordo(\im\,\zeta)^2)\int_{\mathcal{D}\cap\D}
\frac{1}{{(\imag\pi(\zeta-\bar\xi)})^3}\,
\mu(\e^{\imag\pi\xi})\,\diff A(\xi)+\ordo(1),
\label{eq:Ngbasic2}
\end{multline}
as $\im\,\zeta\to0$.
Effectively this reduces to the study of the the derivative of the Bergman projection on 
the upper half-plane $\mathbb{H}$, since the main kernel in \eqref{eq:Ngbasic2} appears
this way. Moreover, in $\mathbb{H}$, we have access to both translation and dilatation 
invariance. This is key to making the argument with boxes effective.

Some clarifying words should perhaps be added regarding the boxes $(1-\epsilon)Q$ compared
with $Q$. More precisely, how should the size reduction be made? It is easier to explain these
matters in the context of the upper half-plane. We work with the hyperbolic metric 
$\diff s_{\Hyp}(\zeta)=(\im \zeta)^{-1}|\diff \zeta|$ and weighted area element 
$\diff A_{-1}(\zeta)=(\im\zeta)^{-1}\diff A(\zeta)$. A prototypical box $Q$ of size $L$ would look like 
this:
\[
Q=\big\{\zeta\in\Hyp: \,0<\re\,\zeta<L,\,\,\e^{-L}<\im\,\zeta<1\big\},
\]
Other boxes could be a horizontal translate of this one, or a dilate of such a translate. The weighted
area of $Q$ is
\[
|Q|_{A_{-1}}=\int_Q\frac{\diff A(\zeta)}{\im\zeta}=\pi^{-1}L^2.
\]
We choose $(1-\epsilon)Q$ to be the smaller box
\[
(1-\epsilon)Q=\big\{\zeta\in\Hyp: \,\epsilon L<\re\,\zeta<(1-\epsilon)L,\,\,
\e^{-(1-\epsilon)L}<\im\,\zeta<\e^{-\epsilon L}\big\}.
\]
Then the weighted area of the box $(1-\epsilon)Q$ is given by
\[
|(1-\epsilon)Q|_{A_{-1}}=\int_{(1-\epsilon)Q}\frac{\diff A(\zeta)}{\im\zeta}=\pi^{-1}(1-2\epsilon)^2L^2.
\]
Moreover, the hyperbolic distance from any point of $(1-\epsilon)Q$ to any point of the complement
$\Hyp\setminus Q$ is at least $\ge\epsilon L$.

\end{document}